\documentclass[11pt,a4paper]{article}

\usepackage[english]{babel}
\usepackage{amsmath,amssymb,amsfonts,a4}
\usepackage{eucal,bbm,stmaryrd}
\usepackage{hyperref}
\usepackage{xcolor}

\begin{document}

\numberwithin{equation}{section}

\newtheorem{Thm}{Theorem}
\newtheorem{Prop}{Proposition}
\newtheorem{Def}{Definition}
\newtheorem{Lem}{Lemma}
\newtheorem{Rem}{Remark}
\newtheorem{Cor}{Corollary}
\newtheorem{Con}{Conjecture}

\newcommand{\Thmautorefname}{Theorem}
\newcommand{\Propautorefname}{Proposition}
\newcommand{\Defautorefname}{Definition}
\newcommand{\Lemautorefname}{Lemma}
\newcommand{\Remautorefname}{Remark}
\newcommand{\Corautorefname}{Corollary}
\newcommand{\Conautorefname}{Conjecture}

\newcommand{\Pf}{\noindent{\bf Proof: }}
\newcommand{\qed}{\hspace*{3em} \hfill{$\square$}}

\newcommand{\N}{\mathbbm{N}}
\newcommand{\Z}{\mathbbm{Z}}

\newcommand{\E}{\mathbbm{E}}
\renewcommand{\P}{\mathbbm{P}}

\newcommand{\Bin}{\mathit{Bin}}

\newcommand{\la}{\lambda}
\newcommand{\La}{\Lambda}
\newcommand{\si}{\sigma}
\newcommand{\al}{\alpha}
\newcommand{\be}{\beta}
\newcommand{\ep}{\epsilon}
\newcommand{\ga}{\gamma}
\newcommand{\de}{\delta}
\newcommand{\De}{\Delta}
\newcommand{\ph}{\varphi}
\newcommand{\om}{\omega}
\renewcommand{\th}{\theta}
\newcommand{\vth}{\vartheta}
\newcommand{\ka}{\kappa}

\newcommand{\pa}{\partial}

\newcommand{\cZ}{\mathcal{Z}}
\newcommand{\cT}{\mathcal{T}}
\newcommand{\cX}{\mathcal{X}}
\newcommand{\cY}{\mathcal{Y}}
\newcommand{\cW}{\mathcal{W}}
\newcommand{\cA}{\mathcal{A}}
\newcommand{\tcW}{\tilde{\cW}}
\newcommand{\cV}{\mathcal{V}}
\newcommand{\cC}{\mathcal{C}}
\newcommand{\cN}{\mathcal{N}}
\newcommand{\da}{\dagger}
\newcommand{\Mda}{M^{\dagger}}
\newcommand{\tM}{\tilde{M}}
\newcommand{\xda}{x_\da}
\newcommand{\xst}{x_*}

\newcommand{\Aut}{Aut}

\newcommand{\Nm}{N_{\wedge}}

\newcommand{\stm}{\setminus}
\newcommand{\lra}{\leftrightarrow}
\newcommand{\lraa}[1]{\stackrel{#1}{\lra}}
\newcommand{\Ra}{\Rightarrow}
\newcommand{\Lra}{\Leftrightarrow}
\newcommand{\bs}{\bigskip}

\makeatletter 
\newcommand{\nobr}{\par\nobreak\@afterheading\vspace{0.3 cm}} 
\makeatother

\thispagestyle{plain}
\title{Bunkbed conjecture for complete bipartite graphs\\ 
and related classes of graphs}
\author{Thomas Richthammer \footnote{Address: Institut f\"ur Mathematik, Universit\"at Paderborn, Warburger Str. 100, 33098 Paderborn, Germany, Email: richth'at'math.uni-paderborn.de
}
}
\maketitle
\begin{abstract}
Let $G = (V,E)$ be a simple finite graph. 
The corresponding bunkbed graph $G^\pm$
consists of two copies  $G^+ = (V^+\!,E^+),G^- = (V^-\!,E^-)$ of $G$ 
and additional edges connecting any two vertices $v_+ \in V_+,v_- \in V_-$ that are the copies of a vertex $v \in V$. 
The bunkbed conjecture states that for independent bond percolation on $G^\pm$, for all $v,w \in V$, 
it is more likely for $v_-,w_-$ to be connected 
than for $v_-,w_+$ to be connected. 
While recently a counterexample for the bunkbed conjecture 
was found, it should still hold for many interesting classes of graphs, and here we give a proof for complete bipartite graphs, complete graphs minus the edges of a complete subgraph, and 
symmetric complete $k$-partite graphs. 
\smallskip 

\end{abstract}





\section{Introduction}

In this introduction we state the bunkbed conjecture 
and review some of the pertinent results. 
We start with a short recap of independent bond percolation. 
Let $G = (V,E)$ be a (simple, undirected) graph, 
and let $p_e \in [0,1]$ for $e \in E$ be a family of edge weights. 
Let $\cZ_e, e \in E$, be independent, $\{0,1\}$-valued random variables (on some probability space) 
with $\P(\cZ_e = 1) = p_e$ for all $e \in E$. 
An edge $e \in E$ is called open if $\cZ_e = 1$ and closed if $\cZ_e = 0$. 
The stochastic process $(\cZ_e)_{e \in E}$ is called 
independent bond percolation on $G$ with edge probabilities $(p_e)_{e \in E}$. 
We may think of percolation as the random subgraph of $G$ obtained from $G$ by removing all closed edges. 
The connected components of this random subgraph are called clusters 
and for $v,w \in V$ we write $v \lra w$, if $v,w$ are contained in the same cluster, i.e. if $v$ and $w$ are joined by a path consisting of finitely many open edges. 
Usually percolation is considered on infinite graphs with some regularity, such as infinite lattices, 
and often the edge weights are assumed to be constant, 
i.e. for some $p \in [0,1]$ we have 
$p_e = p$ for all $e \in E$. We refer to \cite{G} for an introduction to percolation. 

We restrict our attention to percolation on bunkbed graphs. 
For a given graph $G = (V,E)$, the corresponding bunkbed 
graph $G^\pm = (V^\pm,E^\pm)$ is defined by 
$V^\pm := V \times \{+,-\}$ and $E^\pm := \{u_+v_+,u_-v_-: uv \in E\} \cup \{u_+u_-: u \in V\}$, 
where we have written $v_+ := (v,+)$ and $v_- := (v,-)$ and used the 
usual edge notation, e.g. $uv := \{u,v\}$. 
Sometimes $V^+ := V \times \{+\}$ and $V^- := V \times \{-\}$ are called upstairs and downstairs layer respectively, 
and edges of the form $u_+v_+$ or $u_-v_-$ are called horizontal edges, 
whereas edges of the form $u_+u_-$ are called vertical edges.
For given weights $p_G = (p_x)_{x \in E \cup V} \in [0,1]^{E \cup V}$ 
we consider corresponding weights $(p_e)_{e \in E^\pm}$ on $G^\pm$ 
such that 
$p_{u_+v_+} = p_{u_-v_-} := p_{uv}$ and $p_{u_+u_-} := p_u$
for all $u,v \in V$. 
Slightly abusing notation, we also write $p_G$ for $(p_e)_{e \in E^\pm}$. 
We will use $p_G$ as edge probabilities in our percolation model on $G^\pm$, and we write  $\P_{p_G}$ instead of $\P$ 
to indicate the choice of edge weights in our notation of probabilities. 
We note that the reflection map 
$\ph_\pm: V^\pm \to V^\pm$, $\ph_\pm(u_+) := u_-$, 
$\ph_\pm(u_-) := u_+$ for all $u \in V$, 
is a graph automorphism of $G^\pm$, 
which preserves weights $p_G$ as defined above. 
Thus our independent percolation model on $G^\pm$ 
is invariant under this reflection. 

Two special cases of the above model are of interest:  
Given $p \in [0,1]$ and $H \subset V$, we may consider 
probability weights given by $p_{uv} := p$ for all $uv \in E$, 
$p_u := 1$ for all $u \in H$ and $p_u := 0$ for all $u \in V \setminus H$, i.e. we have constant weights on horizontal edges and deterministic vertical edges. 
Alternatively, given $p \in [0,1]$, we may consider 
constant probability weights given by $p_{uv} := p$ for all 
$uv \in E$ and $p_u := p$ for all $u \in V$. 
If we consider these special cases we write $\P_{p,H}$ or  $\P_p$ respectively. 
For the above models we consider the following properties: 

\begin{Def} Bunkbed properties. 
Let $G = (V,E)$ be a (simple, undirected) finite graph. 
We consider the following properties of $G$ 
w.r.t. independent percolation on the bunkbed graph $G^\pm$:  
\begin{itemize}
\item[(B1)] $\forall \; p_G  \in [0,1]^{E \cup V}, v,w \in V: \qquad\P_{p_G}(v_- \! \lra w_-) \ge \P_{p_G}(v_- \! \lra w_+)$. 
\item[(B2)] $\forall \; p \in [0,1], H \subset V, v,w \in V: \quad 
\P_{p,H}(v_- \! \lra w_-) \ge \P_{p,H}(v_- \!\lra  w_+)$. 
\item[(B3)] $\forall \; p \in [0,1], v,w \in V: \qquad \qquad\, \P_{p}(v_- \! \lra w_-) \ge \P_{p}(v_- \! \lra w_+)$. 
\end{itemize}
\end{Def}
In each case the two connection probabilities do not change if 
$v,w$ are interchanged and/or $+,-$ are interchanged, 
which easily follows from the symmetry of $\lra$ and 
the reflection symmetry of the probability models.  
The above properties thus compare the connection probability 
for two vertices in the same layer with that of the two corresponding vertices in different layers. 
The properties are motivated by the heuristics that vertices in different layers are - in some sense - further apart from each other than the corresponding vertices in the same layer, and increasing the 'distance' between vertices in general should make it harder for the vertices to be connected w.r.t. independent percolation. 
\begin{Rem} Bunkbed properties. \nobr
\begin{itemize}
\item 
The above bunkbed properties are obviously related: 
Any graph satisfying (B1) also satisfies (B2), 
and any graph satisfying (B2) also satisfies (B3). 
The latter is even true for fixed $p$, which 
follows via conditioning on all vertical edges. 
For other relations between the properties, see \cite{RS}. 
\item 
Horizontal edges $e$ with $p_e = 0$ may be removed and 
horizontal edges with $p_e = 1$ may be contracted, so one 
may want to exclude these values in (B1). 
Similarly the inequalities in (B2) and (B3) trivially hold in case of $p \in \{0,1\}$. 
\item 
Properties (B1),(B2) and (B3) can also be considered for infinite 
graphs. However, restricting our attention to finite graphs is no loss of generality, since connection probabilities in infinite graphs 
can be obtained as limits of connection probabilities in 
suitable sequences of finite graphs. 
\end{itemize}
\end{Rem}

By the above heuristics, intuitively one might expect
that every finite graph satisfies (B1), (B2) and (B3), 
and this was first conjectured by Kasteleyn (1985) 
as remarked in \cite{BK}. 
In this context (B2) and (B3) are sometimes referred to 
as the strong and weak version of the bunkbed conjecture. 
While the bunkbed conjecture (in any version) was believed 
to hold for a long time, recently it has been disproven:  
\cite{GPZ} shows that (B2) and (B3) do not hold for 
carefully constructed graphs, building on earlier related counterexamples in \cite{H}. 

In light of the recent counterexamples, positive results might be considered to be even more valuable, but only few of these could be obtained so far: 
It is easy to see that a version of (B2) 
holds for $p$ sufficiently small (depending on $G,H$), 
and similarly a version of (B2) holds for $p$ sufficiently large (depending on $G,H$), see \cite{HNK}. 
There are also a few results for special classes of graphs.  
For graphs with very few connections, an inductive approach can be used - e.g. see the strategies used in \cite{L1}, \cite{L2}.
It is not too difficult to see that (B1) holds for all trees and cycle graphs. 
Furthermore, (B2) has been shown for complete graphs, see \cite{HL} 
and earlier partial results in \cite{B}. 
The result for complete graphs mainly relies on 
the high degree of symmetry of the these graphs. 
Our aim is to deal with situations that are somewhat less symmetric. 
We will present our results in the following section. The sections after that are devoted to the proofs of our results.

\section{Results}

Let us first review the result \cite{HL} of Hintum and Lammers for the complete graph. 
It is useful to note that the proof of this result does not use 
the full symmetry of the complete graph. Indeed, we have the following:

\begin{Thm} \label{Thm:old} 
The case of neighboring vertices with a local symmetry. 
Let $G = (V,E)$ be a finite complete graph 
with weights $p_G \in [0,1)^E \times [0,1]^V$. 
Let $v,w \in V$ such that $p_{vw} > 0$. 
In case of $p_v,p_w \neq 1$ suppose that for all 
$u \in V$ with $p_{uw} > 0$ and $p_{u} \neq 1$ 
there is an automorphism $\ph$ of the weighted graph $G$ 
such that $\ph(vw) = uw$. 
Then we have 
$$
\P_{p_G}(v_- \lra w_-) \ge \P_{p_G}(v_- \lra w_+).
$$
\end{Thm}
For the convenience of the reader we have included the proof of Theorem \ref{Thm:old} in Section 5, 
but we note that this proof follows the argument given in \cite{HL} (even if our presentation makes it look somewhat different). 
\begin{Rem} The case of neighboring vertices with a local symmetry. \nobr 
\begin{itemize}
\item 
In our formulation of Theorem 1 we have excluded
$p_e = 1$ in order to simplify the local symmetry condition. 
We note that for any edge with $p_e = 1$ its two 
incident vertices can be combined into a single vertex without changing 
connectivity probabilities, so this is not really a restriction. 
\item 
The local symmetry is formulated in terms 
of the existence of automorphisms. 
We note that an automorphism of a graph $G = (V,E)$ with probability weights $p_G$ 
is simply a graph automorphism $\ph: V \to V$ that 
preserves the probability weights in that 
$p_{\ph(v)} = p_v$ for all $v \in V$ and 
$p_{\ph(u)\ph(v)} = p_{uv}$ for all $uv \in E$. 
What is actually needed in the proof is the weaker 
(but less explicit) condition that the probabilities 
$\P_{p_G}(u_-\lra w_-)$ and  $\P_{p_G}(u_-\lra w_+)$
are constant for $u \in V$ with $p_{uw} > 0$ and $p_u \neq 1$. 
\item 
The theorem is formulated for the model $\P_{p_G}$, and it immediately implies corresponding results for $\P_{p,H}$ and $\P_p$. 
In order to see how the assumptions simplify, let us formulate the 
result for $\P_p$:  Let $G = (V,E)$ be a finite graph and $p \in [0,1]$. 
Let $v,w \in  V$ such that $vw \in E$. 
Suppose that for all $u \in V$ such that $uw \in E$   
there is an graph automorphism $\ph$ of $G$ such that
$\ph(vw) = uw$. 
Then we have 
$\P_p(v_- \lra w_-) \ge \P_p(v_- \lra w_+).$
\end{itemize}
\end{Rem}
The scope of Theorem \ref{Thm:old} is somewhat limited 
by the assumption that $v,w$ are neighbors. 
However, the theorem readily implies that every complete graph has property (B2), and in case of constant probability weights
one can even go beyond this graph class: 
Theorem  \ref{Thm:old} easily implies: 
\begin{Cor} 
Let $G = (V,E)$ be a finite, edge-transitive graph. 
Then we have 
$$
\forall p \in [0,1], vw \in E: \P_p(v_- \lra w_-) \ge \P_p( v_- \lra w_+). 
$$
\end{Cor}
Most interesting examples of edge-transitive graphs 
are symmetric (i.e. vertex- and edge-transitive) or indeed arc-transitive. For relations between the different notions of graph transitivity and examples, we refer to \cite{GR}. 
In order to give some concrete examples, we note that in particular the above corollary applies to the following interesting 
classes of graphs: complete graphs, cycle graphs, 
the graphs corresponding to regular polyhedra (such as the icosahedral graph), hypercube graphs (such as the cubic graph), 
the Petersen graph. It also applies to complete bipartite graphs. 
Here in fact all pairs of vertices that are non-neighboring have the same set of neighbors. 
This is our motivation for investigating 
the following situation, which is the main result of our paper: 
\begin{Thm} \label{Thm:new} 
The case of vertices with the same neighbors. 
Let $G = (V,E)$ be a finite complete graph with probability 
weights $p_G \in [0,1]^{E \cup V}$. 
Let $v,w \in V$ such that for all $u \in V \stm \{v,w\}$ 
we have $p_{vu} = p_{wu}$. 
Then we have 
$$
\P_{p_G}(v_- \lra w_-) \ge \P_{p_G}(v_- \lra w_+).
$$
\end{Thm}
The proof of Theorem \ref{Thm:new} is presented 
in the next section.
It relies on a suitable decomposition of the events, 
conditioning and symmetrization that allows us 
to use the symmetry assumption. While these ideas are rather 
elementary, the way they have to be combined to give the result is somewhat tricky. 
\begin{Rem}
The case of vertices with the same neighbors. \nobr 
\begin{itemize}
\item 
Here the symmetry assumption can also be formulated in terms of 
graph automorphisms: We assume that the transposition map 
$\ph_{v,w}: V \to V$, 
$
\ph_{v,w}(v) := w, 
\ph_{v,w}(w) := v, \ph_{v,w}(u) := u \text{ for all $u \in V \stm \{v,w\}$}
$ 
is an automorphism of the weighted graph $G$. 
\item 
While Theorem \ref{Thm:old} only applies to vertices with graph distance 1, 
the above result only applies to vertices with graph 
distance $\le 2$. 
Unfortunately the symmetry assumption is again rather strong. 
\item 
Again the theorem immediately implies corresponding results for $\P_{p,H}$ and $\P_p$. 
E.g. let us formulate the 
result for $\P_p$:  
Let $G = (V,E)$ be a finite graph and $p \in [0,1]$. 
Let $v,w \in V$ such that for all $u \in V \stm \{v,w\}$ 
we have $vu \in E$ iff $wu \in E$. Then we have  
$\P_p(v_- \lra w_-) \ge \P_p(v_- \lra w_+)$.
\end{itemize}
\end{Rem}
%
%
The combination of Theorems \ref{Thm:old} and \ref{Thm:new} 
proves the bunkbed conjecture for new classes of graphs. 

\begin{Cor} \label{Cor:bipartite}
Let $G = (V,E)$ be a finite complete graph, let $V_1,V_2$ be a disjoint composition of $V$, let $H \subset V$, 
and let $p,p' \in [0,1)$. 
Define probability weights $p_G \in [0,1]^{E \cup V}$ 
by $p_{uv} = 0$ for all $u \neq v \in V_1$, $p_{uv} = p$ for all $u \in V_1,v\in V_2$.  
$p_{uv} = p'$ for all $u \neq v \in V_2$, $p_u = 1$ for all $u \in H$
and $p_u = 0$ for all $u \notin H$. 
Then we have 
$$
\forall v,w \in V: \P_{p_G}(v_- \lra w_-) \ge \P_{p_G}(v_- \lra w_+). 
$$
In particular the bunkbed property (B2) holds for all complete bipartite graphs, 
and for all complete graphs minus the edges of an arbitrary complete subgraph. 
\end{Cor}
In fact the corollary still holds for arbitrary weights $p_v \in [0,1]$, $v \in V$ instead of the deterministic weights above, which easily 
follows by conditioning on the states of all vertical edges. 
We have not found other interesting classes of graphs 
for which a combination of the two theorems implies 
the bunkbed property (B2). For the weaker bunkbed property (B3) 
we can go a little further: 
\begin{Cor} \label{Cor:kpartite} 
Let $k \ge 1$ and let $G = (V,E)$ be a symmetric complete $k$-partite graph, i.e. suppose that 
$V$ can be decomposed into disjoint subsets $V_1,...,V_k$ such that $|V_1| = ... = |V_k|$ 
and $E = \{uv: \exists i \neq j: u \in V_i, v \in V_j\}$. 
Then $G$ satisfies the bunkbed property (B3). 
\end{Cor}
The straightforward proofs for the two preceding corollaries 
are relegated to Section 4. 

Finally we note that the above theorems prove inequalities as in (B2) or (B3) for many graphs and suitable fixed choices of $H \subset V, v,w \in V$ or $v,w \in V$ respectively. However, 
in order to treat further classes of graphs, something is missing. E.g. for general complete $k$-partite graphs with 
vertex set decomposition $V_1,...,V_k$ the case of $v,w \in V_i$
for some $i$ can be treated as above, but in case of $
v \in V_i, w \in V_j$ for $i \neq j$ the symmetry assumption 
in Theorem \ref{Thm:old} is too strong and does not hold. 
However one might hope that the ideas of the proof of Theorem \ref{Thm:new} could lead to further progress in showing bunkbed properties for more general classes of graphs.

\section{Vertices with the same neighbors: 
Theorem \ref{Thm:new}}

For the proof of Theorem \ref{Thm:new} let $v,w \in V$ 
satisfy the given symmetry assumption. 
We will write $\P := \P_{p_G}$. 
By reflection invariance we have  $\P(v_-\lra w_+) = \P(v_+ \lra w_-)$ and 
 $\P(v_- \lra w_-) = \P(v_+ \lra w_+)$, so we need to prove 
$$
d := \P(v_+ \lra w_+) + \P(v_- \lra w_-) - \P(v_+ \lra w_-) - \P(v_- \lra w_+) \ge 0. 
$$
First we decompose
\begin{align*}
\P(v_+ \lra w_+) &= \P(v_- \lra v_+ \lra w_+ \lra  w_-) +   \P(v_- \not \lra v_+ \lra w_+ \lra  w_-)\\ 
&+  \P(v_- \lra v_+ \lra w_+ \not \lra  w_-) +  \P(v_- \not \lra v_+ \lra w_+ \not \lra  w_-) 
\end{align*}
and similarly for the other three events. 
In the corresponding decomposition of $d$, 
the probability with all four vertices connected 
appears twice with a positive and twice with a negative sign, 
and each probability with three of the vertices connected appears once with a positive and once with a negative sign. Thus 
\begin{align*}
d = \P(v_- \not \lra v_+ \lra w_+ \not \lra  w_-) + \P(v_+ \not \lra v_- \lra w_- \not \lra  w_+)\\ 
-\P(v_- \not \lra v_+ \lra w_- \not \lra  w_+) - \P(v_+ \not \lra v_- \lra w_+ \not \lra  w_-).  
\end{align*}
Next we condition on $\cZ_v$. 
We note that on $\cZ_v = 1$ we have $v_+ \lra v_-$ and thus  
$\P(v_- \not \lra v_+ \lra w_+ \not \lra  w_-|\cZ_{v} = 1) = 0$
and similarly for the other three probabilities. 
Thus it suffices to consider the conditional probabilities w.r.t.
$\P(.|\cZ_{v} = 0)$, 
i.e. w.l.o.g. we may assume $p_v = 0$. Similarly we may assume $p_w = 0$. 
Next we consider the edge $vw$. 
Let $A := \{\cZ_{v_+w_+}=\cZ_{v_-w_-} = 0\}$. 
On $A^c$  we either have $v_+ \lra w_+$ or $v_- \lra w_-$, so 
$$
\P(v_- \not \lra v_+ \lra w_- \not \lra  w_+|A^c) = 0 = \P(v_+ \not \lra v_- \lra w_+ \not \lra  w_-|A^c). 
$$
Thus it suffices to consider the conditional probabilities w.r.t.  $\P(.|A)$, 
i.e. w.l.o.g. we may assume $p_{vw} = 0$.
Next let  $W = \{v_-,v_+,w_-,w_+\}$ and let 
$G^\pm_{vw} = (V^\pm_{vw},E^\pm_{vw})$ denote the subgraph 
of  $G^\pm$ induced by $V^\pm_{vw} := V_+ \cup V_- \stm W$. 
For any decomposition $C = \{C_i: i \in I\}$ of $V^\pm_{vw}$ into disjoint connected sets, 
we let $A_C$ denote the event that $\cZ_e, e \in E_{vw}^\pm,$ 
produce clusters given by $C$. 
Let $\P_C := \P(.|A_C)$ denote the conditional distribution. 
Let 
\begin{align*}
d_C = \P_C(v_- \not \lra v_+ \lra w_+ \not \lra  w_-) + \P_C(v_+ \not \lra v_- \lra w_- \not \lra  w_+)\\ 
-\P_C(v_- \not \lra v_+ \lra w_- \not \lra  w_+) - \P_C(v_+ \not \lra v_- \lra w_+ \not \lra  w_-).  
\end{align*}
It suffices to prove $d_C \ge 0$ for all $C$, and we fix $C$ for the remainder of the proof. 
We now look at these four events from the perspective of the clusters $C_i$. 
For $u \in W$ and $i \in I$ 
we write  $u \sim C_i$ if $\cZ_{uu'} = 1$ for some $u' \in C_i$. 
$$
\cA_i = \{u \in W: u \sim C_i\} 
$$
denotes the (random) set of vertices of $W$ directly connected to $C_i$. We note that $v_- \not \lra v_+ \lra w_+ \not \lra  w_-$ 
iff for all $i$ we have $|\cA_i| \le 1$ or $\cA_i = \{v_-,w_-\}$ or $\cA_i = \{v_+,w_+\}$, and the latter occurs at least once. 
(Here we use that $p_v = p_w = p_{vw} = 0$, so all connections
between the points of $W$ have to go via one of the clusters $C_i$.) Thus we have 
\begin{align*}
&\P_C(v_- \not \lra v_+ \lra w_+ \not \lra  w_-)\\
&= \!\!\!\!\sum_{J,K,L: K \neq \emptyset} \!\!\!
\P_C(\forall i \in J\!: |\cA_i| \le 1, \forall i \in K\!: 
\cA_i = \{v_+,w_+\}, \forall i \in L\!: \cA_i =  \{v_-,w_-\})\\
&= \!\!\!\sum_{J,K,L: K \neq \emptyset} \; \prod_{i \in J} 
\P_C(|\cA_i| \le 1) \prod_{i \in K} \P_C(\cA_i = \{v_+,w_+\}) \prod_{i \in L} \P_C(\cA_i =  \{v_-,w_-\}), 
\end{align*}
where the sum is over all disjoint decompositions of the index set $I$. We note that in the last step we may use  independence, 
since the edge sets connecting $W$ to $C_i$ 
and the edge set $E^\pm_{vw}$ are disjoint.  
We have similar decompositions 
for the other three terms in $d_C$. With 
$$
p_M(\le 1) := \prod_{i \in M} \P_C(|\cA_i| \le 1) \quad 
\text{ and } \quad 
p_{M}(W') := \prod_{i \in M} \P_C(\cA_i = W'))
$$
for $M \subset I, W' \subset W$, we thus obtain 
\begin{align*}
&d_C = \sum_{J,K,L: K \neq \emptyset} 
p_J(\le 1) d_{K,L}, \quad \text{ where }\\
&d_{K,L} := p_K(v_+,w_+) p_L(v_-,w_-) + p_K(v_-,w_-) p_L(v_+,w_+)\\
&\hspace*{3 cm} - p_K(v_+,w_-) p_L(v_-,w_+)- p_K(v_-,w_+) p_L(v_+,w_-). 
\end{align*}
Thus it suffices to show that $d_{K,L} \ge 0$ for all 
$K,L$. Writing 
$$
p_i(u) := \P_C(u \sim C_i) \text{ for } u \in W
$$
we note that the given symmetry assumption implies that 
$$
p_i(v_+) = 1 - \prod_{u \in C_i \cap V_+} (1-p_{v_+u}) = 
1- \prod_{u \in C_i \cap V_+} (1-p_{w_+u}) = p_i(w_+) =: p_{i+}
$$
and similarly $p_i(v_-) = p_i(w_-) =: p_{i-}$.
Thus we can write 
\begin{align*}
&d_{K,L} =  \prod_{i \in K}
p_{i+}^2(1-p_{i-})^2\prod_{i \in L} p_{i-}^2(1-p_{i+})^2 + 
\prod_{i \in K} p_{i-}^2(1-p_{i+})^2 \prod_{i \in L}
p_{i+}^2(1-p_{i-})^2 \\
&\qquad - 2\prod_{i \in K}  p_{i-}p_{i+} (1-p_{i-})(1-p_{i+})
\prod_{i \in L} p_{i-}p_{i+} (1-p_{i-})(1-p_{i+}) \\
&= \Big(\prod_{i \in K}
p_{i+}(1-p_{i-})\prod_{i \in L} p_{i-}(1-p_{i+}) - 
\prod_{i \in K} p_{i-}(1-p_{i+}) \prod_{i \in L}
p_{i+}(1-p_{i-})\Big)^2, 
\end{align*}
which implies that $d_{K,L} \ge 0$ and thus we have completed the proof. \qed

\section{Special classes of graphs: Corollaries \ref{Cor:bipartite} and \ref{Cor:kpartite}}

For the proof of Corollary \ref{Cor:bipartite} let $v,w \in V$. 
In case of $w \in V_1,v \in V_2$ we note that the assertion is trivial for $p = 0$, and for $p > 0$ we use Theorem \ref{Thm:old}
noting that $p_{vw} > 0$. 
If $p_v,p_w \neq 1$ and $u \in V$ is an arbitrary vertex with $p_{uw} > 0$ and $p_{u} \neq 1$, 
then necessarily $u \in V_2$ and the transposition 
$\ph_{v,u}$ of $v$ and $u$ is an automorphism of the weighted graph $G$. We have thus verified the required symmetry assumption. (The case $w \in V_2,v \in V_1$ is the same.) 

In case of $w,v \in V_1$ we use Theorem \ref{Thm:new}
noting $p_{vu} = p_{wu} = 0$ for all $u \in V_1 \setminus \{v,w\}$
and  $p_{vu} = p_{wu} = p$ for all $u \in V_2$. 
Similarly, in case of $w,v \in V_2$ we use Theorem \ref{Thm:new}
noting $p_{vu} = p_{wu} = p'$ for all $u \in V_2 \setminus \{v,w\}$
and  $p_{vu} = p_{wu} = p$ for all $u \in V_1$. 

For the last part of the corollary we may assume 
w.l.o.g. that $p < 1$. We note that 
any complete bipartite graph can be obtained from the 
graph under consideration by setting $p' := 0$, and 
any complete graph minus the edges of an arbitrary complete 
subgraph can similarly be obtained by setting $p' := p$. 
\qed 

\bigskip

For the proof of Corollary \ref{Cor:kpartite} 
let $p \in [0,1]$ and $v,w \in V$. W.l.o.g. $p \in (0,1)$. 
In case of $w \in V_i,v \in V_j$ for $i \neq j$ we use Theorem \ref{Thm:old}. Noting that the given graph is edge-transitive, 
the symmetry assumption is satisfied. 
In case of $v,w \in V_i$ for some $i$ we use Theorem \ref{Thm:new} noting that for $u \in V$ we have 
$vu \in E \Lra u \notin V_i \Lra wu \in E$. 
\qed

\section{Neighboring vertices: Theorem \ref{Thm:old}}

For the proof of Theorem \ref{Thm:old} let $v,w \in V$ 
satisfy the given symmetry assumption. 
We write $\P = \P_{p_G}$. 
As in the proof of Theorem \ref{Thm:new} 
we may assume that $p_w = 0$ and $p_v \neq 1$. 
We define 
$$
d := 2\sum_{u} c_{uw} (\P(u_- \lra w_-)- \P(u_- \lra w_+)), \text{ where }c_{uw} := - \ln(1-p_{uw}).  
$$
We note that $c_{uw} \in [0,\infty)$ since $p_{uw} \in [0,1)$. 
For all $u$ such that $p_{uw} = 0$ we have $c_{uw} = 0$. 
For all $u$ such that $p_{uw} > 0$ and $p_{u} = 1$ we have 
 $\P(u_- \lra w_-) = \P(u_- \lra w_+)$. 
For all $u$ such that $p_{uw} > 0$ and $p_{u} \neq 1$
we have $\P(v_- \lra w_-)$ $= \P(u_- \lra w_-)$ 
and $\P(v_- \lra w_+) = \P(u_- \lra w_+)$, 
since the distribution of percolation on the bunkbed graph is 
invariant under all automorphisms of the underlying weighted graph. 
Thus 
$$
d = 2 c (\P(v_- \lra w_-)- \P(v_- \lra w_+)), \text{ where }
c := \sum_{u: p_{uw} > 0,p_{u}< 1} c_{uw}.  
$$
Since $c \ge c_{vw}> 0$ it suffices to show that $d \ge 0$. 
As in the proof of Theorem \ref{Thm:new} reflection invariance 
gives 
$$
d = \sum_{u} c_{uw}(\P(u_- \lra w_-) - \P(u_- \lra w_+) + 
\P(u_+ \lra w_+) -  \P(u_+ \lra w_-)). 
$$
Next let $G^\pm_{w} = (V^\pm_{w},E^\pm_{w})$ denote the subgraph 
of  $G^\pm$ induced by $V^\pm_{w} := V_+ \cup V_- \stm \{w_+,w_-\}$. 
For any decomposition $C = \{C_i:i \in I\}$ of $V^\pm_{w}$ into disjoint connected sets, 
let $A_C$ denote the event that $\cZ_e, e \in E_{w}^\pm,$ 
produce clusters given by $C$,  
let $\P_C := \P(.|A_C)$ denote the conditional distribution, 
and let 
\begin{align*}
d_C \! &:= \!
\sum_{u}\! c_{uw}(\P_C(u_- \!\lra \! w_-)  - \P_C(u_- \!\lra \! w_+) + \P_C(u_+ \!\lra \! w_+) - \P_C(u_+ \!\lra \! w_-))\\ 
&= \!
\sum_{u}\! c_{uw}( \P_C(u_- \lra w_- \not \lra w_+)  
 - \P_C(u_- \lra w_+ \not \lra w_-)\\
&\qquad + \P_C(u_+ \lra w_+ \not \lra w_-) - \P_C(u_+ \lra w_- \not \lra w_+)), 
\end{align*}
where we have canceled  $\P_C(v_+ \lra w_+ \lra w_-)$ and 
$\P_C(v_- \lra w_+ \lra w_-)$. 
It suffices to show that $d_C \ge 0$ for all $C$, and we fix 
$C$ for the remainder of the proof. 
For $i \in I$ and  $w' \in \{w_+,w_-\}$ we write $w' \sim C_i$
iff $\cZ_{w'v'} = 1$ for some $v' \in C_i$, and we let 
$$
p_{i+} := \P_C(w_+ \sim C_i)\quad \text{ and } \quad  
p_{i-} := \P_C(w_- \sim C_i).
$$
We note that in case of $u_- \in C_i$ we have 
\begin{align*}
&\P_C(u_- \lra w_- \not \lra w_+) 
= \P_C(w_- \sim C_i, w_+ \not \sim C_i, \forall j \neq i: w_- \not \sim C_j \text{ or } w_+ \not \sim C_j) \\
&= p_{i-}(1-p_{i+}) r_i, \text{ where } r_i := \prod_{j \neq i}
(1- p_{j-}p_{j+}). 
\end{align*}
Here the first step is due to equality of the two events, 
and in the second step we have used independence due to the events 
depending on disjoint edge sets. 
With similar calculations for the other probabilities we obtain 
\begin{align*}
d_C &=\sum_i
\sum_{u_- \in C_i}  c_{uw} (p_{i-}(1-p_{i+}) r_i - p_{i+}(1-p_{i-}) r_i) \\
&\qquad + \sum_i \sum_{u_+ \in C_i}  c_{uw} (p_{i+}(1-p_{i-}) r_i - p_{i-}(1-p_{i+}) r_i))\\
&= \sum_i \Big( (p_{i-}-p_{i+}) r_i \sum_{u_- \in C_i}  c_{uw} 
+ (p_{i+}-p_{i-}) r_i \sum_{u_+ \in C_i}  c_{uw} )
\end{align*}
By definition of the weights $c_{uw}$ we have 
$$
\sum_{u_- \in C_i} c_{uw} =  - \ln \prod_{u_- \in C_i} (1-p_{uw}) =  - \ln (1-p_{i-})
$$
and similarly for the other sum, so 
$$
d_C =  \sum_i r_i (p_{i-} - p_{i+})(-\ln(1-p_{i-}) + \ln(1-p_{i+})), 
$$
which indeed is nonnegative, since $f(x) = - \ln (1-x)$ is increasing. This finishes the proof of the theorem. \qed

\end{document}